\documentclass{gtart_h}

\def\ifplaintex{\expandafter\ifx\csname documentclass\endcsname\relax}

\def\gtp{{\mathsurround=0pt\it $\cal G\mskip-2mu$eometry \&\ 
$\cal T\!\!$opology $\cal P\!$ublications}}  

\def\Addressesr{\bigskip
{\small \parskip 0pt \leftskip 0pt \rightskip 0pt plus 1fil \def\\{\par}
\sl\theaddress\par
\medskip
\rm Email:\stdspace\tt\theemail\hfill\rm Received:\qua\receiveddate \par}}

\def\recd{{\small Received:\qua\receiveddate\ifx\reviseddate\relax
\else\qquad Revised:\qua\reviseddate\fi\par}} 

\def\lognumber#1{\def\thelognumber{#1}}
\def\volumenumber#1{\def\thevolumenumber{#1}}
\def\volumeyear#1{\def\thevolumeyear{#1}}
\def\papernumber#1{\def\thepapernumber{#1}}
\def\pagenumbers#1#2{\def\startpage{#1}\def\finishpage{#2}}
\def\published#1{\def\publishdate{#1}}

\def\received#1{\def\receiveddate{#1}}

\def\accepted#1{\def\accepteddate{#1}}

\long\def\asciiabstract#1{\long\def\theasciiabstract{#1}}

\let\\\par\let\thelognumber\relax\let\thevolumenumber\relax
\let\thepapernumber\relax\let\thevolumeyear\relax\let\startpage\relax
\let\finishpage\relax\let\publishdate\relax\let\receiveddate\relax
\let\reviseddate\relax\let\accepteddate\relax\let\theasciititle\relax
\let\theasciiauthors\relax
\let\theasciiabstract\relax

\let\theasciiemail\relax

\ifplaintex
\font\logobig=cmssbx10 scaled 3836
\font\logomed=cmssbx10 scaled 2557
\else
\font\logobig=cmssbx10 scaled 4200
\font\logomed=cmssbx10 scaled 2800
\fi

\long\def\makeagttitle{   
\count0=\startpage
\agt\hfill      
\hbox to 45truept{\vbox to 0pt{\vglue -13truept{\logomed A\kern -.37em{\logobig 
T}\kern -.38em G}\vss}\hss}
\break
{\small Volume \thevolumenumber\ (\thevolumeyear)
\startpage--\finishpage\nl
Published: \publishdate}

\vglue .25truein

{\parskip=0pt\leftskip 0pt plus
1fil\def\\{\par\smallskip}{\Large\bf\thetitle}\par\medskip} \vglue
0.05truein

{\parskip=0pt\leftskip 0pt plus 1fil\def\\{\par}{\sc\theauthors}
\par\medskip}%
 
\vglue 0.03truein 

{\small\leftskip 25truept\rightskip 25truept{\bf Abstract}\stdspace\theabstract

{\bf AMS Classification}\stdspace\theprimaryclass
\ifx\thesecondaryclass\relax\else; \thesecondaryclass\fi\par
{\bf Keywords}\stdspace \thekeywords\par}\vglue 7truept

}   

\ifplaintex
\font\phead=cmsl9 scaled 950
\font\pnum=cmbx10 scaled 913
\font\pfoot=cmsl9 scaled 950
\headline{\vbox to 0pt{\vskip -4.5mm\line{\small\phead\ifnum
\count0=\startpage ISSN 1472-2739 (on-line) 1472-2747 (printed)
\hfill {\pnum\folio}\else\ifodd\count0\def\\{ }%
\ifx\theshorttitle\relax\thetitle\else\theshorttitle\fi\hfill{\pnum\folio}
\else\def\\{ and }{\pnum\folio}\hfill\ifx\theshortauthors\relax\theauthors
\else\theshortauthors\fi\fi\fi}\vss}}
\footline{\vbox to 0pt{\vglue 0mm\line{\small\pfoot\ifnum\count0=\startpage
\copyright\ \gtp\hfill\else
\agt, Volume \thevolumenumber\ (\thevolumeyear)\hfill\fi}\vss}}
\else
\font=cmsl9 scaled 1050
\font\lnum=cmbx10 
\font=cmsl9 scaled 1050
\makeatletter
\def\@oddhead{{\small\lhead\ifnum\count0=\startpage ISSN 1472-2739 
(on-line) 1472-2747 (printed)\hfill {\lnum\number\count0}\else\ifodd\count0
\def\\{ }\ifx\theshorttitle\relax \thetitle \else\theshorttitle\fi\hfill
{\lnum\number\count0}\else\def\\{ and }{\lnum\number\count0}
\hfill\ifx\theshortauthors\relax 
\theauthors\else\theshortauthors\fi\fi\fi}}\def\@evenhead{\@oddhead}
\def\@oddfoot{\small\lfoot\ifnum\count0=\startpage\copyright\ \gtp\hfill\else
\agt, Volume \thevolumenumber\ (\thevolumeyear)\hfill\fi}
\def\@evenfoot{\@oddfoot}
\makeatother
\fi
\let\maketitlepage\makeagttitle
\let\makeshorttitle\maketitlepage
\let\maketitle\maketitlepage

\newwrite\gtoutfile
\long\gdef\makeheadfile{  
{\def\\{, }\def\s{ }
\immediate\openout\gtoutfile head.xxx
\immediate\write\gtoutfile{Proxy-for: \ifx\theasciiauthors\relax
\theauthors\else\theasciiauthors\fi\s<\ifx\theasciiemail\relax\theemail\else\theasciiemail\fi>}
\immediate\write\gtoutfile{\noexpand\\}
\immediate\write\gtoutfile{Authors: \ifx\theasciiauthors\relax
\theauthors\else\theasciiauthors\fi}
{\def\\{ }\immediate\write\gtoutfile{Title: \ifx\theasciititle\relax
\thetitle\else\theasciititle\fi}}
\immediate\write\gtoutfile{Subj-class: GT or SG, GR etc}
\immediate\write\gtoutfile{MSC-class: \theprimaryclass\ifx\thesecondaryclass\relax\else, \thesecondaryclass\fi}
\immediate\write\gtoutfile{Journal-ref: Algebraic and Geometric Topology \thevolumenumber\s
(\thevolumeyear) \startpage-\finishpage}
\immediate\write\gtoutfile{Comments: Published by Algebraic and
Geometric Topology at}
\immediate\write\gtoutfile{\s\s\s  http://www.maths.warwick.ac.uk/agt/AGTVol\thevolumenumber/agt-\thevolumenumber-\thepapernumber.abs.html}
\immediate\write\gtoutfile{\noexpand\\}
\immediate\write\gtoutfile{}
\ifx\theasciiabstract\relax
\immediate\write\gtoutfile{\theabstract}\else
\immediate\write\gtoutfile{\theasciiabstract}\fi
\immediate\write\gtoutfile{}
\immediate\write\gtoutfile{\noexpand\\}
\immediate\write\gtoutfile{}
\immediate\closeout\gtoutfile}}  

\def\maketitlepage{\makeagttitle\makeheadfile}
\let\makeshorttitle\maketitlepage
\let\maketitle\maketitlepage

\lognumber{6}
\volumenumber{4}
\volumeyear{2004}
\papernumber{6}
\published{18 February 2004}
\pagenumbers{81}{93}
\received{6 February 2003}
\accepted{20 May 2003}

\usepackage{amsmath,amssymb}
\input xyall

\newtheorem{thm}{Theorem}[section]
\newtheorem{lem}[thm]{Lemma}
\newtheorem{cor}[thm]{Corollary}
\newtheorem{prop}[thm]{Proposition}
\theoremstyle{definition}

\newtheorem{exe}[thm]{Example}

\def\bZ{\mathbb{Z}}

\def\bC{\mathbb{C}}

\def\cG{\mathcal{G}}

\def\oG{\overline{\cG}}

\newcommand{\chft}{chain field theory }
\newcommand{\chftn}[1] {ChFT^{#1 + 1}}
\newcommand{\Lines}{\mbox{Lines}}
\newcommand{\chfts}{chain field theories }

\newcommand{\cGX} {\cG_{n+1}X}

\newcommand{\ra}{\rightarrow}

\DeclareMathOperator{\Hom}{Hom}
\DeclareMathOperator{\Hol}{Hol}

\DeclareMathOperator{\Ext}{Ext}
\DeclareMathOperator{\Ker}{Ker}

\begin{document}


\title{A functorial approach to differential characters}
\author{Paul Turner}
\address{School of Mathematical and Computer Sciences \\Heriot-Watt 
University\\ Edinburgh EH14 4AS\\Scotland}
\email{paul@ma.hw.ac.uk}
\keywords{Differential character, Homotopy quantum field theory}
\primaryclass{57R56}\secondaryclass{53C05, 81T15, 58A10}


\begin{abstract}
We describe Cheeger-Simons differential characters  in terms of a
variant of Turaev's homotopy quantum field theories
based on chains in a smooth manifold $X$.
\end{abstract}
\asciiabstract{%
We describe Cheeger-Simons differential characters  in terms of a
variant of Turaev's homotopy quantum field theories
based on chains in a smooth manifold X.}

\maketitle


\section*{Introduction}
Cheeger-Simons differential characters can be thought of
equivalence classes of some ``higher''
version of line bundles-with-connection. In dimension two this can be
taken to mean gerbes-with-connection, as explained in
\cite{Brylinski:LoopSpacesBook}. One way to think about ``higher''
line bundles-with-connection is in terms of Turaev's homotopy quantum
field theories \cite{Turaev:HomotopyFieldTheoryInDimension2}
(see also
\cite{ro:hqft,BrightwellTurner:RepresentationOfHomotopySurfaceCategory}),
where in dimension two such a thing provides a vector bundle over the 
free loop space
together with a generalised (flat) connection where parallel transport is
defined across surfaces. To make contact with gerbes and differential 
characters
one needs to define a more rigid variation of  1+1-dimensional homotopy 
quantum field
theory as explained \cite{BunkeTurnerWillerton:GerbesHQFT} (see also
\cite{Segal:TopologicalStructures} and for a similar approach
\cite{MackaayPicken:HolonomyAndParallelTransport}).

There is, however, an intrinsic difference between differential
characters  and homotopy
quantum field theories. The former are
defined in terms of {\em homological} information and the latter in
terms of {\em bordism}.  In dimension two this difference is
unimportant (cf.\ the isomorphism between degree two homology and
bordism) but in higher dimensions
one would expect this difference to become apparent. The underlying
geometrical picture of homotopy quantum field theories is however very
appealing: one thinks of a bundle over some space of $n$-manifolds in $X$
with a generalised connection where parallel transport is defined
across $(n+1)$-cobordisms. The motivation for the present work was to
reconcile this picture with the homological needs of differential 
characters.

The ``functorial approach'' of the title refers to the fact that many
geometrical constructions can be defined in terms of representations
of a geometrical category i.e.\ functors from a geometrical category to
a category of vector spaces.  Homotopy quantum field theories are a
good example, but more familiar is the case of a line
bundle-with-connection on $X$. One can define the path category of $X$
as the category with objects the points of $X$ and morphisms smooth
paths between points. A line bundle with connection on $X$ can be
thought of as a functor from the path category of $X$ to the category
of one-dimensional vector spaces: a point in $X$ is assigned to its
fibre and a path to parallel transport along that path. This functor
must also be continuous in an appropriate way.  This was the point of
view in \cite{BunkeTurnerWillerton:GerbesHQFT} where the authors gave
similar description for gerbes-with-connection, by considering rank
one representations of a category with objects loops in $X$ and
morphisms equivalence classes of surfaces in $X$.

We recall now the definition of Cheeger-Simons differential characters
\cite{cheegersimons83}. Letting $Z_{n+1}X$ denote the group of smooth
$(n+1)$-cycles in $X$, a degree $n+1$ differential
character is a homomorphism $f\colon Z_{n+1}X \ra U(1)$ together with
a closed $(n+2)$-form $c$ such that if $\beta$ is an $(n+2)$-chain
then
\[
f(\partial \beta) = \exp (2 \pi i \int_\beta c).
\]
The collection of these is denoted $\widehat{H}^{n+1}(X)$ where the
index $n+1$ follows the convention in \cite{cheegersimons83} (rather
than that in \cite{Brylinski:LoopSpacesBook} where the index $n+2$ is used
for this group).

\subsection*{Outline of the paper}
In order to marry the homological nature of differential characters with the
functorial viewpoint we introduce new objects which we have dubbed
{\em chain field theories}. These are symmetric monoidal functors
from a category whose objects are smooth $n$-cycles in $X$  and whose
morphisms are $(n+1)$-chains
in $X$, to one-dimensional vector spaces. Such an object should be
thought of as a line bundle over
the group of $n$-cycles in $X$ together with a generalised connection
in which parallel transport is defined across $(n+1)$-chains. The
holonomy of such a bundle is a Cheeger-Simons differential
character. The reader should beware that the bundle analogy only goes
so far as we do not demand continuous functors (see also the remarks at
the end of section 2). From one point of view, a chain field theory
provides a possible interpretation of an $n$-gerbe-with-connection.

In Section 1  we define the chain category of $X$, give the
definition of chain field theory and give two important examples. In
Section 2 we prove the following theorem.  \setcounter{section}{2}
\setcounter{thm}{1}\addtocounter{thm}{-1}
\begin{thm}
On a finite dimensional smooth manifold there is an isomorphism from
the group of $(n+1)$-dimensional chain field theories (up to
isomorphism) to the group of $(n+1)$-dimensional differential characters.
\end{thm}

In Section 3 we characterise flat chain field theories as those
that are invariant under deformation by $(n+2)$-chains and finally we
discuss the classification of flat theories by the group
$H^{n+1}(X;U(1))$.

\setcounter{section}{0}

\section{Chain Field Theories}\label{sec:chft}

We will
construct a symmetric monoidal category, $\cG_{n+1}X$ of $n$-cycles
and $(n+1)$-chains in $X$, and then define a chain field theory to be
a 1-dimensional representation of this category. Throughout we will
work with cubical chains, for consistency with the work of Cheeger and
Simons.

\subsection*{Chain categories}

Let $X$ be a smooth finite dimensional manifold. Let $C_kX$ denote the
group of smooth $k$-chains in $X$ and let $Z_kX$ (resp. $B_kX$) be the 
subgroup
of smooth cycles (resp. boundaries).

The $(n+1)$-dimensional {\em chain chain category} of $X$, denoted
$\cG_{n+1}X$ is defined in the following way. The objects are smooth
$n$-cycles in $X$ and a morphism from $\gamma$ to $\gamma^\prime$ is a
smooth $(n+1)$-chain $\sigma$ satisfying $\partial \sigma = -\gamma +
\gamma^\prime$. The composition $\sigma \circ \sigma^\prime$ is
defined to be sum of chains $\sigma + \sigma^\prime$. Associativity
follows from the fact that $C_{n+1}X$ is a group. Noting that the
endomorphisms of an object $\gamma$ can be identified with the group
of $(n+1)$-cycles, we take the zero cycle as the identity morphism for
$\gamma$. To simplify notation we will write $\cG$ for $\cG_{n+1}X$
where there is no ambiguity and we will write $\cG(\gamma,
\gamma^\prime)$ for the set of morphisms from $\gamma$ to
$\gamma^\prime$. We will also make no notational distinction between
the identity morphisms for different $n$-cycles.

We define a bifunctor $\otimes\colon \cG
\times \cG \ra \cG$ on objects by $\gamma_1\otimes \gamma_2 =
\gamma_1 + \gamma_2$, where the sum on the right is taken in $Z_{n}X$ 
and on morphisms
by $\sigma_1 \otimes \sigma_2 =
\sigma_1 + \sigma_2$, where the sum is taken in $C_{n+1}X$. This provides
$\cG$ with the structure of a monoidal category where the monoidal unit is
the zero cycle in $Z_{n}X$.

\begin{prop}
$\cG$ is a strict symmetric strict monoidal groupoid. Its connected
components are in one-to-one correspondence with
$H_{n}(X;\bZ)$.
\end{prop}

\begin{proof}
That $\cGX$ is strict symmetric strict monoidal follows easily from the
fact that $C_{n+1}X$  and $Z_nX$ and abelian groups.

To see that $\cG$ is a groupoid, let $\sigma\in
\cG(\gamma,\gamma^\prime)$ and note that  $-\sigma\in \cG
(\gamma^\prime, \gamma)$ since
$ \partial (-\sigma) = -\partial (\sigma) = - (-\gamma + \gamma^\prime) =
-\gamma^\prime + \gamma$.
Moreover $(-\sigma)\circ \sigma = \sigma + (-\sigma) = 0$ which is the
identity element in $\cG(\gamma, \gamma)$.

To prove the statement about connected components observe that
$\gamma$ is in the same path component as $\gamma^\prime$ if and only
if  there
exists an $(n+1)$-chain $\sigma\in \cG(\gamma, \gamma^\prime)$ such that 
$\partial
\sigma = - \gamma + \gamma^\prime$ i.e.\ $\gamma$ and $\gamma^\prime$
are homologous.
\end{proof}

In fact, the objects of this category also possess inverses and $\cGX$
is a categorical group i.e.\ a group object in the category of groupoids.

\subsection*{The definition of chain field theories}

We let $\Lines$ denote the category with objects 1-dimensional complex
vector spaces with Hermitian inner product and morphisms isometries.
We regard this as a monoidal category under tensor product. For
background information on monoidal categories, functors and so forth
we refer to the appendix in \cite{BunkeTurnerWillerton:GerbesHQFT}
where all relevant definitions can be found.

An {\em $(n+1)$-dimensional chain field theory} on $X$ is a symmetric
monoidal functor $E\colon \cGX \ra \Lines$ together with a closed
differential $(n+2)$-form $c$ such that for any $(n+2)$-chain $\beta$ 
the following
holds:
\[
E(\partial \beta)(1) = \exp (2\pi i \int_\beta c).
\]
The left hand side of this equation should be interpreted in the
following manner. The boundary of an $(n+2)$-chain is an $(n+1)$-cycle
and hence a morphism in $\cG (0,0)$. Since $0$ is the monoidal unit in
$\cG$ and the functor $E$ is monoidal there is an isomorphism $E(0)
\cong \bC$, and in this way $E(\partial\beta)$ is a unitary map $\bC
\ra \bC$. This condition should be thought of as a smoothness
condition of the functor $E$. We note that as part of the definition
of a monoidal functor there are natural isomorphisms
$\Phi^E_{\gamma,\gamma^\prime}\colon E(\gamma)\otimes E(\gamma^\prime)
\ra E(\gamma + \gamma^\prime)$ for objects $\gamma$ and
$\gamma^\prime$.

We say that a \chft is  {\em flat} if the
$(n+2)$-form $c$ is zero.
The reader should think of a chain field theory as a line bundle over
the space of $n$-cycles with parallel transport defined across
$(n+1)$-chains. At first sight it is tempting to provide a more
general definition in which the functor $E$ takes values in the category
of hermitian vector spaces (rather than one-dimensional
ones). However, the objects of $\cG$ have inverses and $E$ is monoidal
so for an object $\gamma$ we have $E(\gamma) \otimes E(-\gamma) \cong
E(0) \cong \bC$ from which it follows that $E(\gamma)$ is one
dimensional.

Two chain field theories are isomorphic when there is a monoidal
natural isomorphism between them. Recall that this requires a natural
transformation $\Psi\colon E\ra E^\prime$ such that
for each object $\gamma$, the map $\Psi_\gamma\colon E(\gamma) \ra
E^\prime (\gamma)$ is an isomorphism and for each pair of objects
$\gamma$ and $\gamma^\prime$
\[
\Psi_{\gamma +  \gamma^\prime} \circ \Phi^E_{\gamma,\gamma^\prime}
= \Phi^{E^\prime}_{\gamma,\gamma^\prime} \circ
( \Psi_\gamma \otimes \Psi_{\gamma^\prime}).
\]
The set of isomorphism classes of $(n+1)$-dimensional chain field
theories on $X$ becomes a group, denoted $\chftn n (X)$, with product
$\star$ defined as follows. Given two
theories $E$ and $F$ form $E\star F$ by defining $(E\star
F)(\gamma) = E(\gamma) \otimes  F(\gamma)$ and $(E\star
F)(\sigma) = E(\sigma) \otimes  F(\sigma)$. The $(n+2)$-form of
$E\star F$ is the sum in the group of $(n+2)$-forms and the monoidal
structure isomorphisms are the obvious ones. The identity of the group
is the
trivial chain field theory, which assigns all objects to $\bC$ and all
morphisms to the identity map. The $(n+2)$-form of the trivial chain
field theory is the zero form and the monoidal structure isomorphisms
are the canonical identification of $\bC \otimes \bC$ with $\bC$.
The inverse of $E$ is defined by
$E^{-1}(\gamma) = E(\gamma)^* = \Hom (E(\gamma), \bC)$ and
$E^{-1}(\sigma) = E(\sigma)^*$. The set of flat \chfts forms a
subgroup of this group.

A chain field theory has the following very useful invertibility
property.  Given a morphism $\sigma$ we have $E(-\sigma) = E(\sigma)^{-1}$.
This is because
\[
E(-\sigma) = E(-\sigma)\circ E(\sigma) \circ E(\sigma)^{-1}\!\!= E(-\sigma + 
\sigma) \circ E(\sigma)^{-1}\!\!= E(0)\circ E(\sigma)^{-1}\!\! =E(\sigma)^{-1}\!.
\]
Just as line bundles with connection have holonomy defined for closed
paths, a \chft has holonomy defined for closed $(n+1)$-chains i.e.\
$(n+1)$-cycles. If $\sigma$ is an $(n+1)$-cycle then it can be
regarded as an element of $\cG (0,0)$ and we define the {\em holonomy}
of $\sigma$ by
\[
\Hol^E(\sigma) = E(\sigma)(1).
\]
Notice that flat theories have trivial holonomy on boundaries since if
$\beta$ is an $(n+2)$-chain then $ \Hol^E(\partial \beta) = \exp (2\pi
i \int_\beta c) = 1$.

  If $\sigma\in\cG(\gamma, \gamma)$ then $\partial \sigma = -\gamma +
\gamma = 0$ so we can also regard  $\sigma$ as an element of $\cG
(0,0)$ and hence holonomy can be defined. As the next lemma shows,
this holonomy is consistent with the map $E(\sigma) \colon E(\gamma)
\ra E(\gamma)$.

\begin{lem}
If  $\gamma$ is an object in $\cG$ and $\sigma$ is an automorphism of
$\gamma$, then the map $E(\sigma)\colon E(\gamma) \ra E(\gamma)$ is
given by
multiplication by $\Hol^E(\sigma)$.
\end{lem}

\begin{proof}
Since $E$ is a monoidal functor there is an isomorphism $\Phi\colon
E(\gamma)\otimes E(-\gamma) \cong E(0) = \bC$. By naturality of the
monoidal structure isomorphisms we have the following commutative
diagram.
\[
\xymatrix{
E(\gamma) \otimes E(-\gamma) \ar[r]^(.65){\Phi} \ar[d]_{E(\sigma) \times
Id} & E(0)
  \ar[d]^{E(\sigma)} \\
E(\gamma) \otimes E(-\gamma) \ar[r]_(.65){\Phi} & E(0)
}
\]
Letting $a$ and $b$ be generators of $E(\gamma)$ and $E(-\gamma)$
respectively we can write  $E(\sigma)( a) = \lambda a$. By chasing 
$a\otimes b$
around the diagram one way we get $\Phi (a\otimes b) \Hol^E(\sigma)$
and the other way $\lambda \Phi(a\otimes b)$. It
follows that $\lambda= \Hol^E(\sigma)$.
\end{proof}

This lemma has two corollaries which will be useful later on.

\begin{cor}\label{cor:hol}
If $\sigma_1$ and $\sigma_2\in \cG(\gamma, \gamma^\prime)$ and
$\Hol^E(\sigma_1 - \sigma_2)= 1$ then $E(\sigma_1) =
E(\sigma_2)$.
  \end{cor}

\begin{proof}
We have that $ E(-\sigma_2) \circ E(\sigma_1) = E(-\sigma_2 \circ 
\sigma_1) =
E(\sigma_1 - \sigma_2)$ and using the lemma above we see that this is
multiplication by $\Hol^E(\sigma_1 - \sigma_2)=1$ i.e.\  $E(-\sigma_2)
\circ E(\sigma_1) = Id_{E(\gamma)}$. Thus
\[
E(\sigma_1) = Id_{E(\gamma^\prime)} \circ E(\sigma_1) = E(\sigma_2) 
\circ E(-\sigma_2)
\circ E(\sigma_1) = E(\sigma_2) \circ Id_{E(\gamma)} = E(\sigma_2).
\]
\end{proof}

\begin{cor}\label{cor:holhomclass}
For a flat theory the holonomy
of an $(n+1)$-cycle $\sigma$ depends only on the
homology class $[\sigma]\in H_{n+1}(X,\bZ)$.
\end{cor}

\begin{proof}
Suppose $\sigma^\prime = \sigma + \partial \beta$ for some
$(n+2)$-chain $\beta$. Then
\[
\Hol^E(\sigma^\prime)/\Hol^E(\sigma) =
\Hol^E(- \sigma)\Hol^E(\sigma^\prime) =
\Hol^E(-\sigma + \sigma^\prime)  =
\Hol^E(\partial \beta) = 1.
\]
\end{proof}

\subsection*{Examples}
We now give two of examples of chain field theories.

\exe\label{ex:df}\hbox{}~

In the first example we construct an $(n+1)$-dimensional \chft from
an $(n+1)$-form. We let $\Omega^*(X)$ denote the smooth complex
differential forms on $X$ .  By $\Omega^*(X)_{0,
   \bZ}$ we denote the subspace of closed forms which have
periods in $\bZ$. Recall from \cite{cheegersimons83} that
there is an injection
\begin{equation}\label{eq:inject}
\Omega^{k}(X) \ra \Hom (C_kX, U(1))
\end{equation}
given by sending $\omega\in \Omega^{k}(X)$ to the map $\beta\mapsto
\exp (2\pi i \int_\beta \omega)$.

Let $\omega \in \Omega^{n+1}(X)$ and define a \chft $E^\omega\colon \cGX \ra
\Lines$ as follows. For any object $\gamma$ set $E^\omega(\gamma) = \bC$ and
for a morphism $(n+1)$-chain $\sigma$ define $E^\omega(\sigma) \colon 
\bC \ra
\bC$ to be multiplication by $\exp (2\pi i \int_\sigma \omega )$. The
monoidal structure is the canonical one and the $(n+2)$-form $c$ is
taken to be $d\omega$. Using Stokes theorem we see that for any
$(n+2)$-chain $\beta$
\[
E^\omega(\partial \beta) (1) = \exp (2 \pi i \int_{\partial \beta} \omega) =
\exp (2 \pi i \int_{\beta} c)
\]
as required.

As the differential on $(n+1)$-forms is linear this gives rise to a
homomorphism
\begin{equation}\label{eq:map}
\Omega^{n+1}(X) \ra \chftn n (X).
\end{equation}
Notice that if $\omega$ is closed then the \chft constructed above is
flat. Moreover if two closed $(n+1)$-forms differ by an exact form then the
resulting \chfts are isomorphic. To see this let $\omega =
\omega^\prime + d\theta$ for some $\theta\in \Omega^n(X)$. For an
object $\gamma \in Z_nX$ define $\tau_\gamma\colon \bC = E^\omega
(\gamma) \ra E^{\omega^\prime}(\gamma) = \bC$ to be multiplication by
$\exp (2\pi i \int_{-\gamma}\theta )$. This defines a natural
transformation $\tau\colon E^\omega \ra E^{\omega^\prime}$. Thus
(\ref{eq:map}) becomes  a homomorphism
\begin{equation}\label{eq:formtoflatchft}
H^{n+1}(X,U(1)) \ra \mbox{Flat}\chftn n (X).
\end{equation}

\exe\label{ex:cs}\hbox{}~

Now  we construct a \chft
from a Cheeger-Simons differential character. Recall (\cite{cheegersimons83}
and \cite{Brylinski:LoopSpacesBook}) that the
Cheeger-Simons group  of differential characters is defined by:

\begin{align*}
\widehat{H}^{n+1}(X) = \{ f\in \Hom (Z_{n+1}X,U(1)) \; \mid \; & \exists 
  c \in
\Omega^{n+2}_{0,\bZ}(X) \mbox{ such that } \\
  & \forall \beta \in C_{n+2}X, \;
f(\partial \beta) = \exp (2\pi i \int_\beta c)\}
\end{align*}
This group fits in to the following exact sequences:
\begin{gather}\label{eq:cses1}
0\strut\rightarrow H^{n+1}(X,U(1))
   \rightarrow \widehat H^{n+1}(X)   
\stackrel{c}{\rightarrow} \Omega^{n+2}_{0,\bZ}(X)\\
0\strut\rightarrow \Omega^{n+1}(X)/\Omega^{n+1}(X)_{0, \bZ}
   \rightarrow\widehat H^{n+1}(X)
   \rightarrow H^{n+2}(X,\bZ)\rightarrow 0
\end{gather}

Starting with a differential character $f\colon Z_{n+1}X \ra U(1)$
with $(n+2)$-form $c$ we will define an $(n+1)$-dimensional chain
field theory $E^f\colon \cG \ra \Lines$ as follows.

There is a short exact
sequence
\[
0 \ra Z_{n+1}X \stackrel{\iota}{\rightarrow} C_{n+1}X 
\stackrel{\partial}{\rightarrow} B_{n}X \ra 0
\]
which gives rise to an exact sequence
\begin{equation}\label{eq:es}
0\ra \Hom (B_{n}X, U(1)) \stackrel{\partial^*}{\rightarrow} \Hom 
(C_{n+1}X,U(1)) \stackrel{\iota^*}{\rightarrow} \Hom (Z_{n+1}X,U(1))
\ra 0.
\end{equation}
This sequence is exact on the left since $U(1)$ is divisible and it
follows that $\Ext(Z_{n+1}X,U(1))$ vanishes. 

Using this exact sequence
choose a lift $\tilde f \colon C_{n+1}X\ra U(1)$ of $f$ and for
objects set $E^f(\gamma) = \bC$ and for morphisms define
$E^f(\sigma)\colon \bC \ra \bC$ to be multiplication by $\tilde f
(\sigma)$. The monoidal structure is taken to be the canonical one and
the $(n+2)$-form is taken to be $c$.

That this provides a well defined symmetric monoidal functor follows
from the fact that $C_{n+1}X$ is an abelian group. The
condition on $c$ is also immediate since for any $\beta\in C_{n+2}X$ we have
that $\partial \beta\in Z_{n+1}X$ so
\[
E(\partial \beta)(1) = \tilde f (\partial \beta) = f(\partial \beta) =
\exp (2\pi i \int_\beta c).
\]
A priori this construction depends on the choice of lift of $f$,
however another choice  yields an
isomorphic chain field theory. Moreover, the construction above is
additive.

\begin{prop}\label{prop:construction}
The construction above provides a homomorphism of groups $\widehat
H^{n+1} (X)\ra \chftn n (X)$.
\end{prop}

\begin{proof}
To show the construction is independent of the lift, let $\overline f$
be another lift which gives rise to another chain field theory
$\overline E^f$ and claim that $E^f$ is isomorphic to $\overline
E^f$.

Noting that $\tilde f /
\overline f \in \Ker ( \Hom (C_{n+1}X,U(1)) \ra
\Hom (Z_{n+1}X,U(1)))$ by using the exact sequence (\ref{eq:es}) we can 
regard
$\tilde f / \overline f$ as a homomorphism $B_{n} \ra U(1)$.
There is
an exact sequence
\[
0\ra B_nX \ra Z_nX \ra H_n(X,\bZ)\ra 0
\]
and thus (again since $U(1)$ divisible) an exact sequence
\begin{equation}\label{eq:homexact}
0\ra \Hom (H_n(X,\bZ) ,U(1)) \ra\Hom (Z_nX ,U(1)) \ra\Hom ( B_nX ,U(1)) 
\ra 0.
\end{equation}
Thus we can lift  $\tilde f / \overline f$ to a homomorphism
$h\colon Z_nX \ra U(1)$. We now
define a natural transformation $\tau \colon E^f \ra \overline E^f$ as
follows. For an object $n$-cycle $\gamma$ in $\cGX$ define
$\tau_\gamma\colon  \bC = E^f (\gamma) \ra \overline
E^f(\gamma) = \bC $ to be multiplication by
$ h(\gamma)$.  Note that since $h$ is a
homomorphism $\tau$ satisfies $\tau_{\gamma+\gamma^\prime} =
\tau_\gamma \tau_{\gamma^\prime}$ and $\tau_{-\gamma} =
\tau_\gamma^{-1}$.  To show that $\tau$ is natural we must show that
for any morphism $\sigma$ from $\gamma$ to $\gamma^\prime$ we have
$\tau_\gamma (1) \overline f (\sigma) = \tau_{\gamma^\prime} (1)
\tilde f (\sigma)$. This is true since
\[
\overline f (\sigma) / \tilde f (\sigma) = (\overline f  / \tilde
f)(\partial \sigma) = h(-\gamma + \gamma^\prime) = \tau_{-\gamma + 
\gamma^\prime} (1) =
\tau_{-\gamma}(1) \tau_{\gamma^\prime}(1).
\]
Thus, up to isomorphism, the construction above is independent of the
choice of lift.

Finally, to see that we have a homomorphism we must show that for
differential characters $f$ and $f^\prime$ we have an isomorphism $
E^{f+f^\prime} \cong E^f \star E^{f^\prime}$. This follows
immediately from the definition of $\star$ and the fact that if we
have lifts $\tilde f$ and $ \tilde f^\prime$ of $f$ and $f^\prime$ we
can choose the lift of $f+f^\prime$ to be $\tilde f + \tilde
f^\prime$, from which we see that the canonical identification of $\bC
\otimes \bC$ with $\bC$ provides an isomorphism from $E^f \star
E^{f^\prime}$ to $E^{f+f^\prime}$.
\end{proof}

If the $(n+2)$-form $c$ above is zero, then the chain field theory
constructed above is flat and using exact sequence
(\ref{eq:cses1}), we can regard the
differential character as an element of $H^n(X,U(1))$ and there is a
homomorphism
\[
H^{n+1}(X,U(1)) \ra \mbox{Flat}\chftn n (X).
\]
This is the same homomorphism as (\ref{eq:formtoflatchft}). In fact
Example \ref{ex:df} is a special case of Example \ref{ex:cs}, using the fact
that an $n+1$-form $\omega$ determines a differential character by $f=
\exp (2\pi i \int \omega)$ and $c=d\omega$.

It is interesting to compare the example above with the constructions
found in the integration theory of  Freed and Quinn
(\cite{FreedQuinn, Freed}) in the context of
Chern-Simons theory.

\section{Classification by Cheeger-Simons groups}\label{sec:classification}
We now show that equivalence classes of $(n+1)$-dimensional \chfts are
classified by the Cheeger-Simons group $\widehat H ^{n+1} (X)$.  Taking
holonomy of a \chft defines a function
\[
\Hol \colon \chftn n (X) \ra \widehat H ^{n+1} (X).
\]
Recall the notation used before: the holonomy of $E$ is denoted
$\Hol^E$. Using this notation the function $\Hol$ above sends $(E,c)$
to $(\Hol^E,c)$ and this function is a homomorphism of groups since
\[
\Hol^{E\star F} (\sigma)= (E\star F)(\sigma)(1) 
=E(\sigma)(1)F(\sigma)(1) = \Hol^{E}
(\sigma)\Hol^{F} (\sigma).
\]
The proof of the following theorem is a reformulation of the proof of
the main theorem in \cite{BunkeTurnerWillerton:GerbesHQFT}.

\begin{thm}\label{thm:main}
On a finite dimensional smooth manifold there is an isomorphism from
the group of $(n+1)$-dimensional chain field theories (up to
isomorphism) to the group of $(n+1)$-dimensional differential
characters.
\end{thm}

\begin{proof}
   We will show that the holonomy homomorphism $\Hol$ is an isomorphism
   with inverse provided by the homomorphism in Proposition
   \ref{prop:construction}.

   Firstly, we will show that
   $\Ker(\Hol)$ is trivial. Let $E\in \Ker(\Hol)$, so $\Hol^E(\sigma) =
   1$ for all $\sigma\in Z_{n+1} X$. Writing $\underline H$ for the
   category with objects the elements of $H_{n}(X,\bZ)$ only identity
   morphisms we can assign to $E$ a symmetric monoidal functor
   $\underline H \ra \Lines$ as follows. Given objects $\gamma$ and
   $\gamma^\prime$ in the same connected component of $\cG$ there is a
   canonical identification of $E(\gamma)$ with $E(\gamma^\prime)$,
   since if $\sigma_1$ and $\sigma_2$ are both morphisms from $\gamma$
   to $\gamma^\prime$ then $\Hol^E(\sigma_1 - \sigma_2) = 1$ and hence
   by Corollary \ref{cor:hol} $E(\sigma_1) = E(\sigma_2)$. It follows
   from the fact that the connected components of $\cG$ are in
   one-to-one correspondence with $H_{n}(X,\bZ)$ that we can associate
   a line $L_x$ to each $x\in H_{n}(X,\bZ)$. Since the morphisms in
   $\underline H$ are identities, this defines a functor $\underline H
   \ra \Lines$. By choosing representatives for each $x\in
   H_{n}(X,\bZ)$, we can use the monoidal structure isomorphisms of $E$
   to define natural isomorphisms $\Phi_{x,x^\prime} \colon L_x\otimes
   L_{x^\prime} \ra L_{x+x^\prime}$ showing that the functor
   $\underline H \ra \Lines$ is monoidal and moreover symmetric.

   Conversely given a symmetric monoidal functor $\phi \colon \underline
   H \ra \Lines$ we can construct a \chft with trivial holonomy by
   setting $E(\gamma) = \phi ([\gamma])$ and $E(\sigma) = Id$. This 
provides an
   identification of $\Ker(\Hol)$ with the group of symmetric monoidal
   functors $\underline H \ra \Lines$. Using Lemma 6.2 of
   \cite{BunkeTurnerWillerton:GerbesHQFT} reformulated for $U(1)$
   rather than $\bC^\times$, the latter can be identified
   with $\Ext(H_{n}(X,\bZ), U(1))$, but this group is trivial since 
$U(1)$ is
   divisible. We have thus shown that the holonomy homomorphism is
   injective.

To see that $\Hol$ is surjective (and that the homomorphism in
Proposition \ref{prop:construction} provides an inverse) let $f$ be a
differential character and  claim that
$\Hol^{E^f} = f$, where $E^f $ is the chain field theory produced in
Proposition \ref{prop:construction}. This is immediate however, since 
for $\sigma\in
Z_{n+1}X$ we know that $E^f(\sigma)\colon \bC \ra \bC$ is
multiplication by $f(\sigma)$, so as an element
of $U(1)$ we have $\Hol^{E^f}(\sigma) = E^f(\sigma)(1) = f(\sigma)$.
\end{proof}

It is important to note that the theorem above relates {\em
   equivalence} classes of chain field theories with differential
characters.  If, for example, one chooses to interpret 1-dimensional
characters as classifying equivalence classes of line
bundles-with-connection then there is only  an identification of line
bundles-with-connection with chain field theories after quotienting up
to equivalence. One could modify the definition of chain field theory
so that the functor is continuous which would get closer to a
genuinely geometric interpretation, but we haven't done that here. I
am grateful to Simon Willerton and Mark Brightwell for clarifying this
point.

\section{Flat theories}
In this section we show that flat \chfts are characterised by invariance
under deformation by $(n+2)$-chains. This is analogous to the fact
that for
flat line bundles parallel transport is
invariant under deformation by homotopy.

Let $\sigma_1,\sigma_2\in \cG(\gamma, \gamma^\prime)$ and suppose
$\beta$ is an $(n+2)$-chain such that $\partial \beta = -\sigma_1 +
\sigma_2$. We say that a \chft $E$ is {\em invariant under chain 
deformation} if for all such
$\sigma_1, \sigma_2$ and $\beta$ we have $E(\sigma_1) = E(\sigma_2)$.

\begin{prop}
A \chft $E$ is flat if and only if it is invariant under chain deformation.
\end{prop}
\begin{proof}
We remarked after the definition of holonomy that if $E$ is flat then
holonomy is trivial on boundaries. Thus
\[
\Hol^E( -\sigma_1 +
\sigma_2) = \Hol^E(\partial \beta) = 1.
\]
So by Corollary \ref{cor:hol}, we see that $E(\sigma_1) = E(\sigma_2)$ and
hence $E$ is invariant under chain deformation.

Conversely, if $E$ is invariant under chain deformation we
claim that $c=0$. Letting  $\beta\in
C_{n+1}X$ we can write  $\partial \beta = -0 + \partial \beta$ and so the
definition of invariance under chain deformation implies
$E(\partial \beta) = E(0) = Id$. Thus for all $\beta\in C_{n+2} X$
\[
\exp (2\pi i \int_\beta c) = E(\partial \beta)(1) = 1.
\]
Using the injectivity of (\ref{eq:inject}) we conclude that $c=0$.
\end{proof}

This can be rephrased as follows. Define $\oG_{n+1}(X)$ to be the quotient
category obtained from $\cG_{n+1}(X)$ by imposing the following relation
on morphisms. Let $\sigma_1,\sigma_2\in \cG (\gamma, \gamma^\prime)$,
then the relation is
\[
\sigma_1 \sim \sigma_2 \mbox{ iff there exists an
$(n+2)$-chain $\beta$ such that $\partial \beta = -\sigma_1 +
\sigma_2$.}
\]
  Composition is still well defined and the
category inherits a monoidal structure from $\cG_{n+1}(X)$.

The above proposition states that a flat \chft is one that factors
through $\oG_{n+1}(X)$. Moreover it is clear that given a symmetric
monoidal functor $\oG_{n+1}(X) \ra \Lines$ the composite $\cG_{n+1}(X)
\ra \oG_{n+1}(X) \ra \Lines$ is a flat \chft and this assignment is
one-to-one. Hence we have the following theorem.

\begin{thm}
There is a one-to-one correspondence between flat \chfts and symmetric
monoidal functors $\oG_{n+1}(X) \ra \Lines$.
\end{thm}

Corollary \ref{cor:holhomclass} states that the holonomy of a flat
chain field theory factors
through $H_{n+1}(X,\bZ)$ and thus may be thought of as a homomorphism
$H_{n+1}(X,\bZ) \ra U(1)$. One may proceed as in the last section to 
study the
function
\[
\Hol \colon \mbox{Flat } \chftn n (X) \ra \Hom(H_{n+1}(X,\bZ),
U(1)) \cong  H^{n+1}(X,U(1))
\]
to establish that this is an isomorphism of groups with the homomorphism
(\ref{eq:formtoflatchft}) providing an inverse. As the proof is
merely a reformulation of the proof of Theorem \ref{thm:main} and the result
is expected once one knows that theorem (compare with the exact
sequence (\ref{eq:cses1})), we omit the details.

\section*{Acknowledgement}
The author thanks Mark Brightwell, Marco Mackaay, Roger Picken  and
Simon Willerton for helpful comments.

\Addressesr

\end{document}